\documentclass[12pt]{amsart}
\usepackage{amssymb, amsmath,amsthm}

\newtheorem{theorem}{Theorem}
\newtheorem{cor}{Corollary}
\newtheorem{lem}{Lemma}
\newcommand{\eps}{\varepsilon}

\begin{document}
\parindent0 in
\parskip 1 em
\title{Discrete spectrum of quantum tubes}
\author{Christopher Lin and Zhiqin Lu}
\date{Dec. 18, 2005}

\email[Christopher Lin, Department of Mathematics, UC Irvine, Irvine, CA 92697]{clin@math.uci.edu}
\email[Zhiqin Lu, Department of Mathematics, UC Irvine, Irvine, CA 92697]{zlu@math.uci.edu}

\thanks{
The
second author is partially supported by  NSF Career award DMS-0347033 and the
Alfred P. Sloan Research Fellowship.}

\maketitle

  \quad A 
quantum tube is essentially a tubular neighborhood about an immersed 
complete manifold in some Euclidean space.  To be more precise, let 
$\varSigma \hookrightarrow \mathbb{R}^{n+k}$, $k\geq 1$, 
$n = dim(\varSigma)$, be an isometric immersion, where $\varSigma$ is 
a complete, noncompact, orientable manifold.  
Then consider the resulting normal bundle 
$T^{\bot}\varSigma$ over $\varSigma$, and the submanifold 
$F = \{(x,\xi) | x\in \varSigma, |\xi < r\} \subset T^{\bot}\varSigma$ for $r$ 
small enough.  The quantum tube is defined as the Riemannian manifold 
$(F, f^*(ds_{E}^2))$, where $ds_{E}^2$ is the Euclidean metric in 
$\mathbb{R}^{n+k}$ and the map $f$ is defined by $f(x,\xi) = x + \xi$.  
If $k=1$, then the quantum tube is also called the quantum layer.
The 
immersion of $\varSigma$ means that the resulting image of $F$ under $f$ in 
$\mathbb{R}^{n+k}$ can have intersections.  Moreover, since $\varSigma$ can 
have quite complicated topology in general, $f(F)$ can too.  However, by doing 
our analysis on $F$ directly (with the pull-back metric), these complications 
are naturally bypassed (cf.~\cite{CEK-1, LL-1}).  

\quad Although on noncompact, noncomplete manifolds there 
is no unique self-adjoint extension of the Laplacian acting on 
compactly supported functions, we can 
always, via the Dirichlet quadratic form define the \text{$\mathit{Dirchlet \hskip 0.1cm Laplacian}$} 
$\Delta_{D}$, which is the self-adjoint extension that reduces to the 
self-adjoint Laplacians defined on complete manifolds and compact manifolds with 
Dirichlet boundary conditions.  Therefore we can proceed to perform spectral 
analysis, in particular, on the quantum tube.  Geometers, like physicists, 
are first and foremost interested in the existence and distribution 
of the discrete  spectrum.  For noncompact manifolds this is in general not 
an easy task at all.  However, using standard variational techniques, the 
authors Duclos, Exner, and Krej\v{c}i\v{r}\'{\i}k were able to, in an interesting 
paper \cite{dek}, prove the existence of discrete spectra for the 
quantum layer (corresponding to $n=2$ and $k=1$ in our defintion) under certain 
integral-curvature conditions on $\varSigma$.  Since the discrete spectrum 
are isolated eigenvalues of finite 
multiplicity, their result is even better, especially in the physical sense since 
the discrete spectrum is composed of energy levels of bound states of a 
nonrelativistic particle.  Our definition of the quantum tube  improved 
theirs in \cite{dek} and we were able to generalize the same existence result 
to the quantum tube.  The challenges  in our attempt at generalization were 
mainly geometrical,  as we sought to replace the necessary geometric conditions with 
appropriate higher dimensional analogs so that similar variational techniques 
from \cite{dek} can be applied meaningfully.  One notable observation that arised 
is the sharp contrast between parabolic and non-parabolic manifolds.  

The main result in~\cite{LL-1} is  as follows:

\begin{theorem}\label{main1}
Let $n\geq 2$ be a natural  number. Suppose $\Sigma\subset
R^{n+1}$ is a complete immersed parabolic  hypersurface such that 
the second fundamental form $A\rightarrow 0$ at infinity.
Moreover, we assume that
\begin{equation}\label{133}
\sum_{k=1}^{[ n/2]}\mu_{2k}{\rm Tr}(\mathcal R^k)\,\, 
\text{is integrable and }\quad
\int_\Sigma \sum_{k=1}^{[n/2
]}\mu_{2k}{\rm Tr}(\mathcal R^k)
      d\varSigma \leq  0,
\end{equation}
where 
$\mu_{2k}>0$ for $k\geq 1$ are positive computable coefficients; $[ n/2 ] $ is the integer part of $n/2$,
and $\mathcal R^k$ is the induced endomorphism of
$\Lambda^{2k}(T_x\Sigma)$ by the curvature tensor $\mathcal R$
of $\Sigma$.
Let $a$ be a positive real number such that $a||A||<C_0<1$
for a constant $C_0$.
If $\Sigma$ is
not totally geodesic, then the ground state of the quantum
layer $\Omega$ exists.
\end{theorem}

In~\cite{LL-2}, we generalized the above results to high codimensional cases:
\begin{theorem}\label{thm3}
Let $(F, f^*(ds_{E}^{2}))$ be an order-$k$ quantum tube with radius r 
and base
manifold $\varSigma$ of dimension $n$ such that the second fundamental form goes to zero at infinity. Moreover, 
we assume that $\varSigma$ is a parabolic manifold, $\sum_{p=1}^{[n/2]} \mu_{2p}
{\rm Tr}(\mathcal R^k)$ integrable, 
 and 
\begin{equation}\label{star}
\int_{\varSigma} \sum_{p=1}^{[n/2]} \mu_{2p}
{\rm Tr}(\mathcal R^k)d\Sigma \, \leq \, 0, 
\end{equation}
If $\varSigma$ is not totally geodesic, then the ground state of the quantum tube from $\Sigma$ exists.
\end{theorem}

By applying the above result into two  dimensional case, we get

\begin{cor}\label{main2}
Suppose that $\Sigma$ is a complete immersed surface of
$R^{n+1}$ such that the second fundamental form $A\rightarrow
0$. Suppose that the Gauss curvature is integrable and suppose
that 
\begin{equation}
e(\Sigma)-\sum\lambda_i\leq 0,
\end{equation}
where $e(\Sigma)$ is the Euler characteristic number of $\Sigma$; $\lambda_i$ is the isoperimetric constant
at each end of $\Sigma$, defined as 
\[
\lambda_i=\underset{r\rightarrow\infty}{\lim}\,\frac{{\rm vol}(B(r))}{\pi r^2}
\]
at each end $E_i$.
Let $a$ be a positive number such that $a||A||<C_0<1$.
If $\Sigma$ is not totally geodesic, then the ground
state of the quantum layer $\Omega$ exists. In particular, if
$e(\Sigma)\leq 0$, then the ground state exists.
\end{cor}

  We  remark here that in the proof of 
Theorem \ref{thm3} (and so as in Theorem \ref{main1} and the analogous result in 
\cite{dek}), the asymptotically flat condition on $\varSigma$ ensures that we get 
a lower bound on the bottom of the essential spectrum, while condition \ref{star} 
(along with parabolicity) enabled us to show that such a bound is also 
a strict upper bound for the total spectrum.  In this way, we were able to conclude 
that the discrete spectrum must be non-empty.  It seems intuitive 
that the asymptotically 
flat condition on $\varSigma$ is essential for there to be discrete spectra, since 
only the ``relatively-curved part of $\varSigma$'' located in the ``interior'' of $\varSigma$ will trap a particle.  If $\varSigma$ is curved more-or-less the same 
everywhere, then a particle may be equally likely to be anywhere since the 
``terrain'' is more-or-less indistinguishable everywhere.  The preceding is of 
course a physical intuition coming from the interpretation of our problem as a 
problem in non-relativistic quantum mechanics, however, it serves to motivate 
the idea that other global curvature assumptions similar to (\ref{star}) may also 
provide the existence of ground state on quantum tubes.

From  Corollary~\ref{main2} (and the result in ~\cite{dek}), it is natural to make the following

\newtheorem*{con}{Conjecture}
\begin{con}
Suppose $\Sigma$ is an embedded asymptotically flat surface in $R^3$ which is not totally geodesic and the Gauss curvature is integrable. Then the ground state of the quantum layer
built from $\Sigma$ exists. 
\end{con}

We have partial results in this direction~\cite{Lu}:

\begin{theorem}[Lu] Suppose $\Sigma$ is asymptotically flat but not totally geodesic in $R^3$. If the Gauss curvature of $\Sigma$ is positive, then the ground state exists for the quantum layer.
\end{theorem}

In general, we have the following result:

\begin{theorem}[Lu] Suppose $\Sigma$ is asymptotically flat but not totally geodesic in $R^3$ and suppose the Gauss curvature is integrable. Let $H$ be the mean curvature. If there is an $\varepsilon>0$ such that
\begin{equation}\label{bbb}
\underset{{r\rightarrow\infty}}{\overline{\lim}}\,\frac 1r\left|
\int_{B(r)}Hd\Sigma\right|\geq\varepsilon,
\end{equation}
then the ground state of the quantum layer exists.
\end{theorem}

\quad
Let's make some remarks on the above results. By the work of~\cite{dek}, we only need to prove the conjecture under the assumption that
\[
\int_\Sigma Kd\Sigma>0.
\]
By a result of Hartman~\cite{hart}, we know that
\begin{equation}\label{aaa}
\frac{1}{2\pi}\int_\Sigma Kd\Sigma=e(\Sigma)-\sum\lambda_i.
\end{equation}
Thus we have $e(\Sigma)>0$, or $e(\Sigma)\geq 1$. Let $g(\Sigma)$ be the genus of $\Sigma$, we then know $g(\Sigma)=0$ and $\Sigma$ must be differmorphic to $\mathbb R^2$, which is a very strong topological restriction.

On the other hand, we have the following lemma:

\begin{lem} Under the assumption that $\int_\Sigma Kd\Sigma>0$, there is an $\eps>0$ such that
\[
\underset{{r\rightarrow\infty}}{\overline{\lim}}\,\frac 1r
\int_{B(r)}|H|d\Sigma\geq\varepsilon.
\]
\end{lem}

{\bf Proof.}  Since $\Sigma$ is differmorphic to $\mathbb R^2$, by ~\eqref{aaa}
\[
0<\int_\Sigma Kd\Sigma\leq 2\pi<4\pi.
\]
Thus by a theorem of White~\cite{white}, we get the conclusion.

\qed

We believe ~\eqref{bbb} is true under the same assumption as in the Lemma.

The above results confirmed the belief that the spectrum of the quantum 
tube  only depends on the geometry of $\varSigma$, its base 
manifold.  With regard to the geometry of $\varSigma$ (or any complete, 
noncompact manifold for that matter), the volume growth (of geodesic balls) 
is an important geometric property.  Roughly speaking, complete, noncompact 
manifolds can be separated into those with at most quadratic volume growth and 
those with faster volume growth.  They are termed (very roughly) parabolic and non-parabolic, 
respectively.  It is the property of parabolicity assumed on $\varSigma$ that 
allowed us to prove the existence of discrete spectra on quantum tubes.  
However, if one looks at the hypothesis of Theorem \ref{thm3}, where 
$\varSigma$ is required to have vanishing curvature at infinity while being 
immersed in Euclidean space, it is highly likely that $\varSigma$ will not 
be of at most quadratic volume growth if $\dim(\varSigma) > 2$, hence unlikely 
to be parabolic.  However, one can be sure that the set of base manifolds 
satisfying the hypothesis of Theorem \ref{thm3} is not empty, due to an 
example provided in \cite{LL-1}.  Nevertheless, it is clear that if one were 
to maintain the assumption of asymptotic flatness of $\varSigma$, then one should  
begin paying attention to the situation when $\varSigma$ is non-parabolic.

\quad
Although we do not yet have a result specifically for quantum tubes over 
non-parabolic manifolds, there is the following preliminary result for general (possibly non-parabolic) 
base manifolds (see \cite{LL-3}):
\begin{theorem}\label{thm4}
Suppose $\varSigma$ is not totally geodesic, 
satisfies the volume growth $V(r) \leq Cr^m$, and whose 
second fundamental form $\vec{A}$ goes to zero at infinity and decays like 
$r^{2}\|\vec{A}\| \to 0$ as $r\to\infty$.  Moreover, suppose 
\begin{equation}\label{new}
\lim_{R\to\infty}
\frac{1}{R^{m-2}}\int_{B(R)}\sum_{p=1}^{[n/2]}\mu_{2p}K_{2p} 
\end{equation}
exists(possibly $-\infty$) and strictly less than $-\frac{1}{4}CC_{1}m^2e^2$, 
where $C_1$ is an explicit constant that depends on the dimension of $\varSigma$, 
radius of the quantum tube, and the upper bound on the curvature of $\varSigma$.  
Then the discrete spectrum of the quantum tube with base manifold 
$\varSigma$ is non-empty. 
\end{theorem}

The result above is certainly an overkill if $\varSigma$ is parabolic.  Thus 
we should think of applying it only to the case of non-parabolic $\varSigma$, 
where $m > 2$.  The direct application of the volume growth hypothesis allows 
one to use polynomially decaying test functions to obtain the condition 
on (\ref{new}), and in turn obtain the upper-bound for 
the bottom of the total spectrum. 

\quad
Theorem \ref{thm4} is only a first step towards generalizing the phenomenon of 
localization (we mean this to be the existence of ground state) to quantum 
tubes over non-parabolic manifolds with similar non-positivy assumptions on 
curvature as the parabolic case.  One clearly cites the technical assumption 
on the decay rate of the second fundamental form, and one would like to remove 
it.  In addition, the negativity condition on (\ref{new}) is very strong.  We 
do not yet know if weaker assumptions such as (\ref{star}) are applicable to the 
case where $\varSigma$ is non-parabolic.  

\quad
{\bf Acknowledgement.} This short note is based on the talk the second author given at {\it Hayama Symposium on
Complex Analysis in Several Variables 2005} on December 18-21, 2005. He thanks the organizers, expecially Professor T. Ohsawa for the invitation and  hospitality during the symposium.

\end{document}